\documentclass[graybox]{svmult}

\usepackage{type1cm}        
\usepackage{makeidx}         
\usepackage{graphicx}        
\usepackage{multicol}        
\usepackage[bottom]{footmisc}

\usepackage{newtxtext}       %
\usepackage[varvw]{newtxmath}       

\usepackage{stackrel}

\makeindex             

\begin{document}

\title*{
 About the structural stability of Maxwell fluids: 
 convergence toward elastodynamics
}
\author{S\'ebastien Boyaval}
\authorrunning{Structural stability of Maxwell fluids toward elastodynamics}
\institute{S\'ebastien Boyaval \at LHSV, Ecole des Ponts, EDF R\&D, Chatou, France, \email{sebastien.boyaval@enpc.fr}
\at MATHERIALS, Inria, Paris, France}
%
%
\maketitle

\newcommand{\bA}{\boldsymbol{A}}
\newcommand{\be}{\boldsymbol{e}}
\newcommand{\boldf}{\boldsymbol{f}}
\newcommand{\bF}{\boldsymbol{F}}
\newcommand{\bI}{\boldsymbol{I}}
\newcommand{\bsigma}{\boldsymbol{\sigma}}
\newcommand{\btau}{\boldsymbol{\tau}}
\newcommand{\bx}{\boldsymbol{x}}
\newcommand{\bzero}{\boldsymbol{0}}
\def\bD{\boldsymbol{D}}
\newcommand{\bu}{\boldsymbol{u}}
\newcommand{\RR}{\mathbb{R}} \newcommand{\R}{\mathbb{R}} 
\renewcommand{\div}{\operatorname{div}}
\newcommand{\grad}{\boldsymbol{\nabla}}

\abstract{
Maxwell's 
models for viscoelastic flows are famous 
for their potential 
to unify 
elastic motions of solids with viscous motions of liquids
in the continuum mechanics perspective.
But rigorous proofs are lacking.
The present note is a 
contribution toward well-defined 
viscoelastic flows proved to encompass both solid and (liquid) fluid regimes.
In a first part, we consider the structural stability of \emph{particular viscoelastic flows}: 
1D shear waves solutions to damped wave equations.
We show the convergence toward purely elastic 1D shear waves solutions to standard wave equations,
as the relaxation time $\lambda$ and the viscosity $\dot{\mu}$ grow unboundedly $\lambda\equiv\frac1G\dot{\mu}\to\infty$
in Maxwell's constitutive equation
$$
\lambda \stackrel{\Diamond}{\btau} 
+ \btau = 2 \dot{\mu}  \bD(\bu)
$$
for the stress $\btau$ of viscoelastic fluids with velocity $\bu$.
In a second part, we consider the structural stability of general \emph{multi-dimensional viscoelastic flows}. 
To that aim, we embed Maxwell's constitutive equation in a symmetric-hyperbolic system of 
PDEs which we proposed in our previous publication [ESAIM:M2AN 55 (2021) 807-831] so as to define multi-dimensional viscoelastic flows unequivocally. 
Next, we show the continuous dependence of multi-dimensional viscoelastic flows on $\lambda\equiv\frac1G\dot{\mu}$
using the relative-entropy tool developped for symmetric-hyperbolic systems after C.~M.~Dafermos.
It implies convergence of the viscoelastic flows defined in [ESAIM:M2AN 55 (2021) 807-831]
toward compressible neo-Hookean elastodynamics when $\lambda\to\infty$.
}

\section{
Maxwell fluids as links between solids and Newtonian fluids}
\label{sec:intro}

We consider 
the viscoelastic motions of a 
fluid body occupying on times $t\in[0,T)$
a subset $\Omega \subset \RR^3$ of the Euclidean ambiant space equipped with a Cartesian system of coordinates $\{x^i, i=1\dots 3\}$.

Denoting $\bu=u^i\be_i$ the velocity field of the fluid, $\rho$ the mass density (a scalar field), $\boldf$ a bulk force field, we assume the following mass and momentum balances
\begin{align}
\label{eq:determinant_spatial_tens} 
& \partial_t \rho + \div \left( \rho \bu \right) = 0 
\\
\label{eq:momentum_spatial_tens}
& \partial_t \left( \rho \bu \right) +  \div \left( \rho \bu\otimes\bu - \bsigma \right) = \rho \boldf
\end{align}
using a model of Maxwell type \cite{maxwell-1874} for the extra-stress $\btau$ in Cauchy 2-tensor 
$\bsigma=-p{\boldsymbol\delta}+\btau$,
$p(\rho)$ being a pressure in the fluid and $\boldsymbol\delta$ the identity 2-tensor, i.e.:
\begin{equation}
\label{eq:UCMmodified} 
\lambda \stackrel{\Diamond}{\btau} 
+ \btau = 2 \dot{\mu}  \bD(\bu) \,.
\end{equation}
In Maxwell's constitutive equation \eqref{eq:UCMmodified}, 
$\dot\mu>0$ is a viscosity parameter, 
$\lambda>0$ is a relaxation-time parameter,
and $\stackrel{\Diamond}{\btau}$ is an \emph{objective} time-rate operator
see e.g. \cite{Sar-2016-cfma}.

\smallskip

It is widely admitted that Maxwell fluid models (i.e. those using \eqref{eq:UCMmodified}) can 
link 
fluids 
where $\btau \xrightarrow{\lambda\to0} 2 \dot{\mu} \bD(\bu)$ in the Newtonian limit,
denoting $\bD(\bu)=\frac12\left(\grad\bu+\grad\bu^T\right)$ as usual,
with
solids governed by elastodynamics when $\lambda\equiv\frac1G\dot{\mu}\to\infty$.
Besides, in applications, one often considers the Upper-Convected Maxwell (UCM) model, 
with objective time-rate 
$\stackrel{\Diamond}{\btau}$ in \eqref{eq:UCMmodified} 
defined 
by the Upper-Convected (UC) derivative
\begin{equation}
\label{eq:UC}
\stackrel{\triangledown}{\btau} := \partial_t \btau + (\bu\cdot\grad) \btau - (\grad\bu)\btau  - \btau(\grad\bu)^T
\end{equation} 
since the formal limit of \eqref{eq:UCMmodified} when $\lambda\equiv\frac1G\dot{\mu}\to\infty$ 
i.e. $\stackrel{\triangledown}{\btau}=2G 
\bD(\bu)$
is compatible with Hookean elastodynamics, i.e. the case when
$\btau=\frac{\dot{\mu}}{\lambda}\left(\bF\bF^T-\bI\right)$
and $\bF$ is the deformation gradient associated with $\bu$, governed 
by 
\begin{equation}
\label{eq:F}
(\partial_t +\bu\cdot\grad)\bF=(\grad\bu)\bF \,.
\end{equation}
But we are not aware of a rigorous proof of such a \emph{structural stability} result for Maxwell fluid models
that would link general multi-dimensional motions of (Hookean) solids with Newtonian fluids.

\bigskip

One-dimensional (1D) shear waves can be defined unequivocally with
\eqref{eq:determinant_spatial_tens}--\eqref{eq:momentum_spatial_tens}--\eqref{eq:UCMmodified}.
The Newtonian fluid limit $\btau \xrightarrow{\lambda\to0} 2 \dot{\mu} \bD(\bu)$ of such 1D shear waves
is established in \cite{Payne2001} as a consequence of the structural stability of \emph{linear} UCM models.
But we are not aware of a rigorous proof of the elastodynamics limit $\lambda\equiv\frac1G\dot{\mu}\to\infty$, though.

Multi-dimensional time-continuous motions cannot be defined unequivocally with
\eqref{eq:determinant_spatial_tens}--\eqref{eq:momentum_spatial_tens}--\eqref{eq:UCMmodified}
in general, even on a small time interval using smooth initial conditions on the whole space $\RR^3$,
because the quasilinear system \eqref{eq:determinant_spatial_tens}--\eqref{eq:momentum_spatial_tens}--\eqref{eq:UCMmodified} may not be hyperbolic. 
%
So the question of structural stability cannot be properly addressed for general multi-dimensional viscoeastic flows
as such with \eqref{eq:determinant_spatial_tens}--\eqref{eq:momentum_spatial_tens}--\eqref{eq:UCMmodified},
i.e. as one usually ``defines'' viscoeastic flows of Maxwell type in the literature.

\bigskip

In 
Sec.~\ref{sec:1D}, we extend the study of \cite{Payne2001} (for linear Maxwell equations) 
to the convergence toward solid elastodynamics when $\lambda\equiv\frac1G\dot{\mu}\to\infty$, 
specifically for the \emph{1D shear waves} solutions to \eqref{eq:determinant_spatial_tens}--\eqref{eq:momentum_spatial_tens}--\eqref{eq:UCMmodified} which are recalled in Sec.~\ref{sec:waves}. 
Such specific studies are a first step toward 
the 
structural stability of more general (nonlinear) models of Maxwell type,
and to a fully rigorous link between solid and fluid regimes using Maxwell models.

Next, to address the structural stability of physically-relevant (nonlinear) Maxwell models,
a further step is to first unequivocally define \emph{multi-dimensional} motions through solutions to \eqref{eq:determinant_spatial_tens}--\eqref{eq:momentum_spatial_tens}--\eqref{eq:UCMmodified}.
We propose here to build upon our former work \cite{Boyaval2021}, thus to consider the structural stability of
\emph{unequivocal viscoelastic flows} defined as solutions to a quasilinear system of PDEs with a symmetric-hyperbolic reformulation that implies \eqref{eq:UCMmodified}.

In a nutshell, our reformulation of Maxwell flows interprets the extra-stress as $\btau=\rho G(\bF\bA\bF^T-\bI)$,
with a view to extending to \emph{Maxwell fluids with finite parameters $\dot{\mu},\lambda\equiv\frac1G\dot{\mu}>0$}
an elastodynamics system where $\bF$ is the deformation gradient associated with 
$\bu$ and $\rho$.
That is the reason why our reformulation \cite{Boyaval2021} requires
\begin{equation}
\label{eq:deformation_spatial_tens}
\partial_t \left( \rho \bF \right) -  \grad \times \left( \rho \bF^T\times\bu \right) = \bzero
\end{equation}
like in elastodynamics, along with the famous involution termed Piola's identity
\begin{equation}
\label{eq:piola_spatial_tens}
\div( \rho \bF^T ) = \bzero \:. 
\end{equation}
Notice that \eqref{eq:deformation_spatial_tens}--\eqref{eq:piola_spatial_tens}--\eqref{eq:determinant_spatial_tens} together
imply \eqref{eq:F}.
Moreover, we assume 
\begin{equation}
\label{eq:Aeul}
\lambda (\partial_t +\bu\cdot\grad)\bA  + \bA = \bF^{-1}\bF^{-T}
\end{equation}
for the symmetric positive definite 
2-tensor $\bA$.
Then, a constitutive equation of Maxwell-type \eqref{eq:UCMmodified} 
holds (for smooth compressible flows \cite{Bollada2012}),
and solutions to the system 
\eqref{eq:determinant_spatial_tens}--\eqref{eq:momentum_spatial_tens}--\eqref{eq:deformation_spatial_tens}--\eqref{eq:piola_spatial_tens}--\eqref{eq:Aeul} without source term,
where $\btau=\rho G (\bF\bA\bF^T-\bI)$ and $p=-\partial_{\rho^{-1}}e_0$ with $e_0$ convex in $\rho^{-1}$,
additionally satisfy a conservation law for a scalar 
quantity that is convex 
in a conserved variable ${\rm U}(\bu,\bF,\rho,\bA)$ (see \cite{Boyaval2021}):
\begin{equation}
\label{eq:energy}
\eta := \frac\rho2 |\bu|^2 + \rho e_0 + \rho \frac{G}2 \bF\bA:\bF \,.
\end{equation}
That is, using notations of \cite[Chap.~V]{DafermosBook4} and denoting $\xi=\frac1\lambda>0$,
there exists a variable change in a convex domain $\rm U \in \mathcal{O}$
such that our system rewrites 
\footnote{ 
 Involutions are keys here, we refer to the short summary \cite{cras2022Boyaval} for instance.
} 
\begin{equation}
\label{eq:symmetric}
\rm
\partial_t U + \partial_\alpha G_\alpha(U)=\xi \Pi(U)
\end{equation}
with smooth vector fluxes $\rm G_\alpha$, 
and smooth fluxes $\rm Q_\alpha(U)$ exist so that $\rm \eta(U)$ satisfies 
\begin{equation}
\label{eq:entropy}
\rm
\partial_t \eta(U) + \partial_\alpha Q_\alpha(U)=\xi D_U \eta(U)\cdot \Pi(U) \,.
\end{equation}
Consequently, our system \eqref{eq:determinant_spatial_tens}--\eqref{eq:momentum_spatial_tens}--\eqref{eq:deformation_spatial_tens}--\eqref{eq:piola_spatial_tens}--\eqref{eq:Aeul}, equiv.~\eqref{eq:symmetric} after a variable change,
admits a symmetric-hyperbolic formulation, see \cite{Boyaval2021},
and one can define unequivocally time-continuous flows of Maxwell fluids on small time intervals
given general smooth initial conditions. 
%
%
So the question of structural stability can be considered for our reformulation of Maxwell fluids,
in particulat using standard results for symmetric-hyperbolic systems \cite{DafermosBook4}.
Note to that aim that, in the hyperbolicity domain $\mathcal{O} \ni \rm U_1, U_2$, the source term $\rm \Pi(U)$ is such that 
\begin{equation}
\label{eq:C0source}
\rm |\Pi_m(U_1)-\Pi_m(U_2)|\le C_m \|U_1-U_2\|^2
\end{equation}
for each component $\rm m=1\ldots 1+d+d^2+\frac{d(d+1)}2$
of the system \eqref{eq:symmetric}.

In Sec.~\ref{sec:full},
we show the continuous dependence on $\lambda\equiv\frac1G\dot{\mu}$,
of general multi-dimensional viscoelastic flows defined unequivocally following \cite{Boyaval2021},
using the relative-entropy tool developped for symmetric-hyperbolic systems after C.~M.~Dafermos.
It implies the following structural stability result: 
convergence of our 
viscoelastic flows 
toward compressible neo-Hookean elastodynamics when 
$\lambda\equiv\frac1G\dot{\mu}\to\infty$.

\section{Setting of the problem for 1D viscoelastic shear waves} 
\label{sec:waves}

A 
shear wave $\bu=u(t,y)\be_x$, $\btau=\tau^{xy}(t,y)\be_x\otimes\be_y$ solution 
to \eqref{eq:determinant_spatial_tens}--\eqref{eq:momentum_spatial_tens}--\eqref{eq:UCMmodified}
on $\{t\ge 0$, $x^2\equiv y\in\Omega\}:=(y_{\min},y_{\max})\subset\RR$,
can be built unequivocally given initial conditions 
$$
u(t=0,y)=u^0(y) \quad \tau^{xy}(t=0,y)=\tau^0(y)
$$
plus boundary conditions at $y\in\partial\Omega$ when necessary (when $y_{\min},y_{\max}$ finite), 
as briefly recalled below. 
Indeed, \eqref{eq:determinant_spatial_tens}--\eqref{eq:momentum_spatial_tens}--\eqref{eq:UCMmodified} 
reduces to: 
\begin{align}
\label{eq:momentum1D}
\partial_t u = \partial_y \tau^{xy} + f^x \,,
\\
\label{eq:maxwell1D}
\lambda \partial_t \tau^{xy} + \tau^{xy} = \dot{\mu} \partial_y u \,, 
\end{align}
when one assumes $\rho$ constant 
(this is natural for the 1D motion along $\be_x\equiv \be_{x^1}$
of a 2D body 
with a Lagrangian description using material coordinates $a=x-X(t,y)$, $b=y\equiv x^2$ in a Cartesian frame,
such that $u\equiv\partial_t X$, $\tau^{xy}\equiv\partial_y X$, 
as it is the case for our reformulation of Maxwell fluids in \cite{Boyaval2021}).

\smallskip

When $\Omega\equiv\{y>0\}$ for instance, the Stokes first 
problem for (\ref{eq:momentum1D}--\ref{eq:maxwell1D}),
with $u^0\equiv0\equiv\tau^0$, 
$u(t,y=0)=UH(t)$ ($U\in\RR^+_*$, $H$ denoting Heaviside's function),
can be solved \emph{analytically} 
by particular 1D shear waves \cite{PREZIOSI1987239}. 
Next, for those particular 1D shear waves, the Newtonian fluid limit $\lambda\to0$
can be established directly 
using an analytical expression of the solution, see e.g. \cite[(4.3)--(4.4)]{JORDAN2004}.

However, the same structural stability result can be established much more generally,
i.e. for a class of well-defined solutions, 
by an \emph{energy method},
i.e. an analysis using energy estimates satisfed by the solutions.
For instance,
the Newtonian fluid limit of well-defined solutions to (\ref{eq:momentum1D}--\ref{eq:maxwell1D})
is a direct 
consequence of the structural stability established in \cite{Payne2001}, 
for a large class of solutions to \emph{linear} Maxwell models (with limited physical relevance, though).

\smallskip

In the sequel, we use arguments similar to \cite{Payne2001} (i.e. energy estimates satisfied by the solutions)
to analyze the structural stability, when $\lambda\equiv\frac1G\dot{\mu}$
or equivalently when $\xi:=\frac1\lambda\to0$ keeping $G$ fixed, of 
the damped wave system (\ref{eq:momentum1D}--\ref{eq:maxwell1D})
which we rewrite using $\tau:=\tau^{xy}(t,y)$, $f:=f^x$ for simplicity:
\begin{align}
\label{eq:momentum1Dbis}
\partial_t u - \partial_y \tau = f
\\
\label{eq:maxwell1Dbis}
\partial_t \tau - G \partial_y u = -\xi \tau \,. 
\end{align}
Note that contrary to the limit $\lambda\to0$ studied in \cite{Payne2001}, 
the limit $\xi:=\frac1\lambda\to0$ studied here
is ``easier'' because it is a \emph{non-singular} limit: only a \emph{lower-order} term vanishes when $\xi\to0$, 
which does not change the hyperbolic type of (\ref{eq:momentum1Dbis}--\ref{eq:maxwell1Dbis}) in the limit.
A detailed proof is nevertheless useful with a view to extending the structural stability result
to general multi-dimensional viscoelastic flows solutions to a complex quasilinear system of PDEs.
Introducing the variables $w^\pm = \tau \pm \sqrt{G} u$, it is also useful to rewrite (\ref{eq:momentum1Dbis}--\ref{eq:maxwell1Dbis}) as
\begin{equation}
\label{eq:riemann}
\partial_t w^\pm \mp \sqrt{G} \partial_y w^\pm = \pm \sqrt{G}f - \frac\xi2 (w^++w^-) \,. 
\end{equation}

\section{Structural stability of 1D shear waves when $\frac1\lambda\equiv\xi\to0$}
\label{sec:1D}

Given $\xi\ge0$ and an open subset $\Omega:=(y_{\min},y_{\max})\subset\RR$,
time-continuous solutions $u(t,y)$, $\tau(t,y)$ to (\ref{eq:momentum1D}--\ref{eq:maxwell1D})
on $t\ge 0$,
with value in $L^2(\Omega)$, 
are well-defined given 
$$
u(t=0,y)=u^0(y)\in L^2(\Omega)
\quad \tau(t=0,y)=\tau^0(y) \in L^2(\Omega)
\quad f(y) \in L^2\left((0,T)\times\Omega\right) 
$$
when $y\in\Omega\equiv\R$.
We recall that $u(t,y),\tau(t,y)$ in fact take values in $H^1(\Omega)\subset L^2(\Omega)$ then, see \cite{DafermosBook4}.

When $y_{\min},y_{\max}$ are finite,
time-continuous solutions $u(t,y)$, $\tau(t,y)$ to (\ref{eq:momentum1D}--\ref{eq:maxwell1D})
remain well-defined on additionally specifying boundary conditions at $y\in\partial\Omega$,
typically \emph{maximally dissipative} 
given $g \in L^2\left((0,T)\times\partial\Omega\right)$ as follows
\begin{align}
\label{eq:bcl} 
z_l^- := & c_l^u u + c_l^\tau \tau = g_l \text{ with $c_l^u c_l^\tau<0$, $c_l^\tau \neq -\sqrt{G} c_l^u $ at } y_{\min}\,,
\\
z_r^+ := & c_r^u u + c_r^\tau \tau = g_r \text{ with $c_r^u c_r^\tau>0$, $c_r^\tau \neq \sqrt{G} c_r^u $ at } y_{\max}\,.
\label{eq:bcr} 
\end{align}
The latter solutions satisfy the following energy estimate
\begin{multline}
\label{eq:wavebound0}
  \frac{d}{dt} \int_\Omega \frac12 \left( |w^+|^2 + |w^-|^2 \right)
  +\frac{\sqrt{G}}{2c_l^u c_l^\tau} |z_l^+|_{y_{\max}}^2
  +\frac{\sqrt{G}}{2c_r^u c_r^\tau} |z_r^-|^2_{y_{\max}}
  \\
  + \xi \int_\Omega (w^++w^-)^2  
  = \int_\Omega f (w^+-w^-) 
  +\frac{\sqrt{G}}{2c_l^u c_l^\tau} g_l^2 +\frac{\sqrt{G}}{2c_r^u c_r^\tau} g_r^2
\end{multline}
where $z_l^- = c_u u - c_\tau \tau$, $z_r^- = c_u u - c_\tau \tau$,
on multiplying \eqref{eq:momentum1D} by $G u$ and \eqref{eq:maxwell1D} by $\tau$.

Consider now two solutions $(u_1,\tau_1)$ and $(u_2,\tau_2)$ for two parameter values $\xi_1\ge\xi_2\ge0$,
with same initial conditions and source term. On $\Omega$, they satisfy
\begin{align}
\label{eq:momentum1Dter}
\partial_t (u_1-u_2) - \partial_y (\tau_1-\tau_2) = 0 \,,
\\
\label{eq:maxwell1Dter}
\partial_t (\tau_1-\tau_2) - G \partial_y (u_1-u_2) 
=  - \xi_1 (\tau_1-\tau_2) - (\xi_1-\xi_2) \tau_2 \,,
\end{align}
and homogeneous boundary conditions on $\partial\Omega$ i.e. (\ref{eq:bcl}--\ref{eq:bcr}) with $g_l=0=g_r$, hence 
\begin{multline}
\label{eq:wavebound1}
  \frac{d}{dt} \int_\Omega \frac12 \left( |w^+_1-w^+_2|^2 + |w^-_1-w^-_2|^2 \right)
  +\frac{\sqrt{G}}{2c_l^u c_l^\tau} |z_{1,l}^+-z_{2,l}^+|_{y_{\max}}^2
  +\frac{\sqrt{G}}{2c_r^u c_r^\tau} |z_{1,r}^+-z_{2,r}^+|^2_{y_{\max}}
\\
  + \xi_1 \int_\Omega  (\tau_1-\tau_2)^2 
  = - (\xi_1-\xi_2) \int_\Omega \tau_2  (\tau_1-\tau_2)
\end{multline}
with obvious notations,
on multiplying \eqref{eq:momentum1Dter} by $G(u_1-u_2)$ and \eqref{eq:maxwell1Dter} by $(\tau_1-\tau_2)$.

Using Cauchy-Schwarz and Young inequalities with \eqref{eq:wavebound1}, one finally obtains:
\begin{proposition}
\label{prop}
Given two parameter values $\xi_1\ge\xi_2\ge0$,
two 
solutions $(u_1,\tau_1)$ and $(u_2,\tau_2)$ of the damped wave system
(\ref{eq:momentum1Dbis}--\ref{eq:maxwell1Dbis}) with same conditions satisfy
\begin{multline}
\label{eq:wavebound2}
  \frac{d}{dt} \int_\Omega \frac12 \left( |w^+_1-w^+_2|^2 + |w^-_1-w^-_2|^2 \right)
  +\frac{\sqrt{G}}{2c_l^u c_l^\tau} |z_{1,l}^+-z_{2,l}^+|_{y_{\max}}^2
  +\frac{\sqrt{G}}{2c_r^u c_r^\tau} |z_{1,r}^+-z_{2,r}^+|^2_{y_{\max}}
\\
  \frac{\xi_1+\xi_2}2 \int_\Omega  (\tau_2-\tau_1)^2 
  \le \frac{\xi_1-\xi_2}2 \int_\Omega \tau_2^2
\end{multline}
i.e. continuous dependence on the relaxation parameter
$$(u_1,\tau_1)\xrightarrow{\xi_1\to\xi_2}(u_2,\tau_2)$$ 
in $L^2(\Omega)$ for all times $t\ge0$, as well as in $L^2(0,T)$ on the boundary $\partial\Omega$.
\end{proposition}

When $\xi_2=0$ in particular, the latter structural stability result of Prop.~\ref{prop} yields convergence
of $(u_1,\tau_1)$ solution of the damped wave system (\ref{eq:momentum1Dbis}--\ref{eq:maxwell1Dbis})
toward $(u_2,\tau_2)$ solution of a standard wave system that coincides with 1D elastodynamics. 

So convergence toward elastodynamics is quite a simple structural-stability result for 
the 1D viscoelastic shear waves solutions to (linear) Maxwell equations -- in comparison with
a singular limit like convergence toward Newtonian fluids \cite{Payne2001}.

However, that non-singular limit is not easy 
anymore when considering nonlinear equations for general multi-dimensional motions,
defined e.g. through our reformulation in \cite{Boyaval2021}.

\section{Structural stability of 
general Maxwell flows when $\frac1\lambda\equiv\xi\to0$}
\label{sec:full}

To establish convergence toward elastodynamics of general multi-dimensional viscoelastic flows,
we consider our reformulation \cite{Boyaval2021} of Maxwell flows using a symmetric-hyperbolic 
system of conservation laws \eqref{eq:symmetric} with variable $\rm U$,
and we use the standard comparison tool introduced by C.~M.~Dafermos: the \emph{relative entropy} \cite[Chap.~V]{DafermosBook4}.
Precisely, consider two classical solutions $\rm U_1, U_2$ using same conditions but two relaxation parameters $\xi_1, \xi_2$.
It holds
\begin{multline}
\rm
\partial_t \left( \eta(U_1)-\eta(U_2)-D_U\eta(U_2)\cdot(U_1-U_2) \right)
\\ \rm
+ \partial_\alpha \left( Q_\alpha(U_1)-Q_\alpha(U_2)-D_U\eta(U_2)\cdot(G_\alpha(U_1)-G_\alpha(U_2)) \right)
\\ \rm
= 
\xi_1 \left( D_U \eta(U_1) - D_U \eta(U_2) \right)\cdot \Pi(U_1) 
\\ \rm
- \left( \partial_t D_U\eta(U_2) \right)\cdot(U_1-U_2)
-\left( \partial_\alpha D_U\eta(U_2) \right)\cdot(G_\alpha(U_1)-G_\alpha(U_2))
\\ \rm
= \xi_1 \left( D_U \eta(U_1) - D_U \eta(U_2) \right)\cdot \Pi(U_1) 
- \xi_2 (U_1-U_2)\cdot D^2_{UU}\eta(U_2) \cdot \Pi(U_2)
\\ \rm
- \partial_\alpha U_2 \cdot D^2_{UU}\eta(U_2) \cdot \left( G_\alpha(U_1)-G_\alpha(U_2) - D_UG_\alpha(U_2)\cdot(U_1-U_2) \right) 
\label{eq:entropyrel}
\end{multline}
which can be compared to \cite[(5.2.10)]{DafermosBook4}:
our relative-entropy equality holds for two \emph{classical} solutions $\rm U_1, U_2$, 
with an additional source term (first line of RHS in \eqref{eq:entropyrel}).

Then, we suggest to compare two well-defined classical solutions $\rm U_1, U_2$ 
that use the same conditions in the hyperbolicity domain $\mathcal O$ 
on complementing \cite[(5.2.14)]{DafermosBook4} as follows
to take into account the additional source term.
In~\eqref{eq:entropyrel} we use:
\begin{multline}
\rm
  \xi_1 \left( D_U \eta(U_1) - D_U \eta(U_2) \right)\cdot \Pi(U_1) 
- \xi_2 (U_1-U_2)\cdot D^2_{UU}\eta(U_2) \cdot \Pi(U_2)
\\ \rm 
=   \left( D_U \eta(U_1) - D_U \eta(U_2) \right)\cdot (\xi_1 \Pi(U_1) - \xi_2 \Pi(U_2))
+ \xi_2 Z(U_1,U_2)\cdot \Pi(U_2)
\\ \rm
=  \xi_1 \left( D_U \eta(U_1) - D_U \eta(U_2) \right)\cdot ( \Pi(U_1) - \Pi(U_2))
\\ \rm
+ (\xi_2- \xi_1) \left( D_U \eta(U_1) - D_U \eta(U_2) \right)\cdot \Pi(U_2)
+ \xi_2 Z(U_1,U_2)\cdot \Pi(U_2)
\end{multline}
where $Z(U_1,U_2)$ is quadratic in $U_1-U_2$.

We conclude using \eqref{eq:C0source} with
Cauchy-Schwarz 
and Young inequalities for any $t\in(0,T),r>0$,
to complement \cite[(5.2.16)]{DafermosBook4} in our case with 
source terms.
Note that by contrast with the linear case of Sec.~\ref{sec:1D}, 
coercivity does not hold i.e.
$$\rm
c_0 \left| D_U \eta(U_1) - D_U \eta(U_2) \right|^2
+ \left( D_U \eta(U_1) - D_U \eta(U_2) \right)\cdot ( \Pi(U_1) - \Pi(U_2))
$$
cannot be guaranteed non-positive 
for all $\rm U_1,U_2$ whatever $c_0>0$~!
Still, $\forall \sigma\in(0,t)$
\begin{multline}
\label{eq:relative}
\int_{|\bx|\le r+c(t-\sigma)} \left| U_1(t)-U_2(t) \right|^2
\le \int_0^\sigma ds \int_{|\bx|\le r+c(t-s)} d\bx \Bigg(
C (\xi_1-\xi_2)^2 |U_2(s)|^2 
\\ + C' \left( 1+ \xi_1^2 \right) |U_1(s)-U_2(s)|^2 \Bigg) 
\end{multline}
holds given positive constants $c,C,C'$ depending solely on the initial conditions and $\rm U_2, \xi_2$
but not $\xi_1$.
So Gromwall's inequality allows one to conclude about the structural stability of general multidimensional Maxwell flows:
%
%
\begin{proposition}
\label{prop2}
Given two parameter values $\xi_1\ge\xi_2\ge0$ bounded above,
consider two smooth solutions $\rm U_1, U_2$ to \eqref{eq:symmetric} on $[0,T)\times\R^3$
with same initial condition of bounded support.
There exists $\rm C_T(U_2)$ such that
$$
\rm \|U_1(t)-U_2(t)\|_{L^2(\R^3)} \le C_T|\xi_1-\xi_2| \quad \forall t\in(0,T) \,.
$$
\end{proposition}

When $\rm \xi_2=0$ in particular, the structural stability result of Prop.~\ref{prop2} yields convergence
of $\rm U_1$ solution to our formulation in \cite{Boyaval2021} of viscoelastic Maxwell flows
toward $\rm U_2$ solution of an elastodynamics system for compressible hyperelastic materials of Neo-Hooken type. 

\section{Conclusion and Perpsectives}

In this short note, we have completed the structural stability results of \cite{Payne2001}
for global-in-time \emph{linear} Maxwell flows (Proposition~\ref{prop} for 1D shear waves)
in the ``easy'' elastodynamics limit case $\lambda\sim\dot{\mu}\to\infty$, 
as opposed to the singular Newtonian fluid limit $\lambda\to0$ in \cite{Payne2001}.
Note that 
in a particular case with smoother solutions, our elastodynamics solid limit case $\lambda\sim\dot{\mu}\to\infty$
was already covered by a structural stability result from \cite{GurGulec2016}.


Moreover, using our recent reformulation of Maxwell's models \cite{Boyaval2021}
so as to unequivocally define viscoelastic flows as solutions to symmetric-hyperbolic PDEs,
we could extend in Prop.~\ref{prop2} the structural stability result to generic multi-dimensional viscoelastic flows:
as widely believed, multi-dimensional viscoelastic flows of Maxwell type can be \emph{rigorously} linked with
solid elastodynamics (of compressible hyperelastic Neo-Hooken materials)
when the relaxation time grows unboundedly.

The result should still hold in the case of \emph{non-isothermal} Maxwell flows defined similarly in \cite{S0219891622500096}
using a symmetric-hyperbolic reformulation which extends our former work \cite{Boyaval2021} (to non-isothermal Maxwell flows,
as well as to viscoelastic flows with finite-extensibility effects and various objective time-rates).
But a second step in structural stability should then furthermo cover more physical liquid-solid transitions, driven by temperature changes.
And it remains a challenge to 
establish structural stability results 
for multidimensional solutions of physically-relevant Maxwell models 
in the \emph{singular} Newtonian-limit case $\lambda\to0$.

%

\begin{acknowledgement}
The author acknowledges the partial support of the ANR project 15-CE01-0013 SEDIFLO: ``Modelling and simulation of solid transport in rivers''.
\end{acknowledgement}
%


\end{document}